\documentclass[a4paper,12pt,legno]{article}

\usepackage{amsmath}
\usepackage{amsfonts}
\usepackage{amsthm}
\usepackage{mathrsfs}
\usepackage{amssymb}

\newtheorem{thm}{THEOREM}
\newtheorem{lem}{LEMMA}

\newtheorem{cor}{COROLLARY}
\newtheorem{ex}{Example}

\theoremstyle{remark}
\newtheorem{rem}{\textbf{Remark}}

\usepackage{graphicx}
\usepackage[T1]{fontenc}
\usepackage[polish]{babel}
\usepackage[utf8]{inputenc}

\newcommand{\Rset}{\mathbb{R}}
\newcommand{\Cset}{\mathbb{C}}

\newcommand{\al}{\alpha}
\newcommand{\be}{\beta}

\begin{document}

\title{Remarks on compositions of some random integral mappings\footnote{Research funded by Narodowe Centrum Nauki (NCN)
Dec2011/01/B/ST1/01257}}

\author{ Zbigniew J. Jurek (University of Wroc\l aw\footnote{Part of this work was done when Author was visiting
Indiana University, Bloomington, \qquad \qquad \qquad \qquad USA in
Spring 2013.}\,\,)}

\date{Published in Stat.Probab. Letters 137 (2018) 277-282}

\maketitle
\begin{quote} \textbf{Abstract.} The random integral
mappings (some type of functionals of L\'evy processes) are continuous 
homomorphisms between convolution subsemigroups of the semigroup of all infinitely divisible measures.  Compositions of those random integrals (mappings) can be always  expressed as another single random integral mapping. That fact is illustrated by some old and new examples.


\medskip


\medskip
\medskip
\emph{Abbreviated title: Compositions of random integral mappings}

\end{quote}
\maketitle

Addresses:

\medskip
\noindent Institute of Mathematics \\ University of Wroc\l aw  \\
Pl. Grunwaldzki 2/4
\\ 50-384 Wroc\l aw \\ Poland. \\ www.math.uni.wroc.pl/$\sim$zjjurek ;
e-mail: zjjurek@math.uni.wroc.pl

\newpage
For the last few decades random integrals were used to describe some classes of limiting distributions. For example, L\'evy class L (selfdecomposable) distributions or  s-selfecomposable distributions (the class U). In those and other situations one had to identify an integrand, an interval (or a half-line) over each one integrates, a class of integrators (L\'evy processes) and a time change in the process; cf. Jurek (2011) for a review of the research in that area and appropriate references;  cf. Jurek (2014) for examples of some specific mappings.

\noindent Here we discuss a quite general set-up for such random integral mappings and prove, among others,  that they are closed under compositions. It is illustrated by some explicit examples.

\medskip
For an interval $(a,b]$ in the positive half-line, two deterministic
functions $h$ and $r$,  and a L\'evy process $Y_{\nu}(t), t\ge0$,
where $\nu$ is the law of random variable $Y_{\nu}(1)$, we consider
the following mapping
\[
\mathcal{D}^{h,r}_{(a,b]}\ni \nu\longmapsto
I^{h,r}_{(a,b]}(\nu):=\mathcal{L}\big(\int_{(a,b]}h(t)\,dY_{\nu}(r(t))\big),
\ \ (\star)
\]
where $\mathcal{L}$ denotes the probability distribution of the random  integral and $ \mathcal{D}^{h,r}_{(a,b]}$ is the domain of the mapping $I^{h,r}_{(a,b]}$ in $(\ast)$.

Our results and proofs are given  for $\Rset^d$ variables. An infinite dimensional case is discussed  in the Concluding  Remarks.

\medskip
\textbf{1. Main results.}

1 a).  \underline{\it{Compositions of the random integral mappings.}}

\medskip
Let time changes $r(t),\, a<t\le b,$ be either $\rho\{s:
s>t\}$ or $\rho\{s: s\le t\}$ for some positive, possibly infinite,
measure $\rho$ on $[0,\infty)$.
For functions $h_1,...,h_m$ on the intervals $(a_1, b_1],..., (a_m,b_m]$, respectively,
and positive measures $\rho_1,...,\rho_m$, let us define
\begin{multline}
\textbf{h}:= h_1 \otimes...\otimes h_m,\ \ \mbox{(the tensor product of
functions)} \\ \mbox{i. e.} \ \
\textbf{h}(t_1,t_2,...,t_m):=h_1(t_1)\cdot
h_2(t_2)\cdot...\cdot h_m(t_m), \ \mbox{where} \ a_i <t_i  \le b_i; \\
\textbf{(a,b]}:=(a_1,b_1]\times ...\times (a_m,b_m], \
\boldsymbol{\rho}:= \rho_1\times...\times \rho_m, \, \mbox{(product
measure)}.
\end{multline}

\begin{thm}
Let functions $h_i$, measures $\rho_i$ (given by increments of
functions $r_i$) and intervals $(a_i,b_i]$, for $i=1,2,...,m$, be as
above.

If the image $\textbf{h((a,b])}=(c,d]\subset \Rset^+$ and $\nu\in
ID(\Rset^d)$ is from an appropriate domain then we have
\begin{equation}
I_{(a_1, b_1]}^{h_1, \rho_1}(I_{(a_2, b_2]}^{h_2,
\rho_2}(...(I_{(a_m, b_m]}^{h_m,\rho_m}(\nu))... ))= I_{(c,d]}^{t, \,
\textbf{h}\,\,\boldsymbol{\rho}}(\nu),
\end{equation}
where $\textbf{h}\,\boldsymbol{\rho}$ is the image of the product
measure $\boldsymbol{\rho}=\rho_1\times ... \times \rho_m$ under the
mapping $\textbf{h}:=h_1 \otimes...\otimes h_m$.

Random integrals $I^{h_i,\rho_i}_{(a_i,b_i]}, i=1,2,...,m$, commute
on the domain
$\mathcal{D}^{t,\,\textbf{h}\,\,\boldsymbol{\rho}}_{(c,d]}$ of  the mapping $I_{(c,d]}^{t, \,
\textbf{h}\,\,\boldsymbol{\rho}}$.
\end{thm}

Examples illustrating Theorem 1 are given in Section 3 below.

\begin{cor}
If $Z_1, Z_2,..., Z_m$ are stochastically independent variables with their probability distributions $\rho_i$ concentrated on 
intervals $(a_i,b_i]$, respectively then
\[
r(t):=\textbf{h}\,\,\boldsymbol{\rho}(s\le t)=P[h_1(Z_1)\cdot ...
\cdot h_m(Z_m)\le t] 
\]
is the time change in the corresponding composition of random mapings $I^{h_i, r_i}_{(a_i, b_i]}$
where $r_i(t)=\rho_i\{s:s\le t\}$.
\end{cor}
If a random mapping is a composition of other mappings we may infer
some inclusions of their ranges. Namely we have
\begin{cor} If  an equality $I^{h,r}_{(a,b]}=
I^{h_1,r_2}_{(a_1,b_1]}\circ I^{h_2,r_2}_{(a_2,b_2]}$ (a
composition) holds on the domain $\mathcal{D}^{h,r}_{(a,b]}$ then we
have
\[
\mathcal{R}^{h,r}_{(a,b]}\equiv
I^{h,r}_{(a,b]}(\mathcal{D}^{h,r}_{(a,b]})\subset
I^{h_1,r_1}_{(a_1,b_1]}(\mathcal{D}^{h_1,r_1}_{(a_1,b_1]})\cap
I^{h_2,r_2}_{(a_2,b_2]}(\mathcal{D}^{h_2,r_2}_{(a_2,b_2]})=\mathcal{R}^{h_1,r_1}_{(a_1,b_1]}\cap
\mathcal{R}^{h_2,r_2}_{(a_2,b_2]}.
\]
\end{cor}

\medskip
1 b). \underline{\it{ Properties  of the random integral mappings. }}

\begin{thm} (a) Assume that $h(a):=h(a+), r(a):=r(a+)$ exist in $\Rset$.
 Then the mapping
\begin{equation}
 \mathcal{D}^{h,r}_{(a,b]}\ni\nu\to I^{h,r}_{(a,b]}(\nu) \in ID
\end{equation}
is a continuous homomorphism between the corresponding measure convolution
semigroups. 

\noindent (b) For given $s>0$, we have that\ $\nu \in
\mathcal{D}^{h,r}_{(a,b]}$ if and only if $\nu^{\ast
s}\in\mathcal{D}^{h,r}_{(a,b]}$,  and
\begin{equation*} 
I^{h,r}_{(a,b]}(\nu^{\ast s})= (I^{h,r}_{(a,b]}(\nu))^{\ast s} = (I^{h,\,sr}_{(a,b]}(\nu)).
\end{equation*}
(c) For $u >0$ and the dilation operator $T_u :\Rset^d\to \Rset^d $ defined as $T_u (x)= u\,x$, we have that  $\nu \in
\mathcal{D}^{h,r}_{(a,b]}$ if and only if $T_u \nu \in
\mathcal{D}^{h,r}_{(a,b]}$, and
\begin{equation*}
T_u\big(I^{h,r}_{(a,b]}(\nu)  \big) =  I^{h,r}_{(a,b]}(T_u\nu)=I^{u
h,r}_{(a,b]}(\nu).
\end{equation*}
(d) For bounded linear operator A on $\Rset^d$ and $\nu \in
\mathcal{D}^{h,r}_{(a,b]}$ we have that $A\nu\in\mathcal{D}^{h,r}_{(a,b]}$ and  
$
A(I^{h,r}_{(a,b]}(\nu))=I^{h,r}_{(a,b]}(A\nu).
$
\end{thm}

\textbf{2. Proofs.}

First we will recall some basic definitions and facts.

\medskip
2 a). \underline{\it{L\'evy-Khintchine representations.}}

\medskip
Here $ID\equiv ID(\Rset^d)$ stands for the class of all infinitely divisible probability measures $\nu$ on $\Rset^d$. Thus their characteristic functions (Fourier transforms) are of the form (the famous L\'evy-Khintchine formula)
\begin{multline}
\hat{\nu}(y)= e^{\Phi(y)}, \  \mbox{where} \  \ \Phi(y)\equiv \log \hat{\nu}(y)=i<y,z> - \frac{1}{2}<y, Ry> \\  + \int_{\Rset^d\setminus{\{0\}}}\,[\,e^{i<y,x>}-1-i<y,x>1_{\{||x||\le 1\}}(x)\,]\,M(dx) , \  y \in \Rset^d,
\end{multline}
 where $<\cdot , \cdot>$  denotes the scalar product and  the triple: a vector $z \in \Rset^d$,  a covariance operator $R$ of a Gaussian part of $\nu$ and a  L\'evy (spectral) measure $M$ (of Poissonian part),  is uniquely determined by $\nu$. In short, we write:  $\nu=[z,R,M]$ and  \  $\Phi$ is referred to as \emph{ the L\'evy exponent} of $\nu \in ID$;  cf. Meerschaert  and Scheffler (2001).  [For more general case than $\Rset^d$,   cf. Araujo-Gine (1980) or Parthasarathy (1967).]

\medskip
2 b). \underline{\it{Definition of path-wise random integrals.}}

\medskip
For an interval ${(a,b]}$ in a positive half-line, a
real-valued continuous of bound variation function $h$ on $(a,b]$, a
positive non-decreasing right-continuous (or non-increasing left-continuous) 
time change function $r$ on $(a,b]$ and a cadlag L\'evy stochastic processes $(Y_{\nu}(t), 0 \le
t < \infty)$, let us define,  via the formal integration by parts
formula, the following \emph{random integral}
\begin{multline}
\int_{(a,b]} h(t)dY_{\nu}(r(t)):=
\\ h(b)Y_{\nu}(r(b))-h(a)Y_{\nu}(r(a))- \int_{(a,b]}Y_{\nu}(r(t)-)dh(t) \in \Rset^d, 
\end{multline}
and the corresponding \emph{random integral mapping}
\begin{equation}
\qquad \qquad  \mathcal{D}^{h,r}_{(a,b]}\ni \nu\to  I^{h,r}_{(a,b]}(\nu):=\mathcal{L}\big( \int_{(a,b]}
h(t)dY_{\nu}(r(t))\big) \in ID,  
\end{equation} 
where  $\mathcal{D}^{h,r}_{(a,b]}$  is the domain of the mapping  $I^{h,r}_{(a, b]}$, that is a
subset of the class ID  consisting of those measures $\nu$ for which
the integral (5) is well defined.

\medskip

From properties of infinitely divisible measures (distributions) one concludes that  the law of the random integral (5)  is  infinitely divisible one; 
cf. Jurek-Vervaat (1983), Lemma 1.1 (or for a particular case,  cf. Jurek-Mason (1993), Section 3.6)).

\begin{rem}

(a) Since L\'evy processes are semi-martingales the random integral (5)
can be defined as an Ito stochastic integral. However, for our purposes we do not need that generality of stochastic calculus.
\noindent [Term: \emph{random integral} emphasizes that h in (5) is a deterministic function ( not a stochastic  process).]

(b) Integrals over intervals (a,b) or (a,$\infty$) or [a,b] and
others are defined as weak limits of integrals over
intervals (a,b] in (5).  
Thus, the random integral $\int_{(a,\infty)} h(t)dY_{\nu}(r(t))$ is well-defined if and only if the function 
\begin{equation}
\ \Rset^d\ni y \to \int_{(a,\infty)} \Phi(h(t)y)dr(t)\in \Cset  
\end{equation}
is a L\'evy exponent  (a functional of the form (3), above).
\end{rem}

2 c). \underline{\it{ L\'evy exponents. }}

\medskip
If $\nu \in \mathcal{D}^{h,r}_{(a,b]}$ and $I^{h,r}_{(a,b]}(\nu)$
have the L\'evy exponents $\Phi$ and $\Phi^{h,r}_{(a,b]}$,
respectively then, from already mentioned in Lemma 1.1 in Jurek-Vervaat(1983), we get
\begin{equation}
\Phi^{h,r}_{(a,b]}(y)=\int_{(a,b]} \Phi(h(t)y)dr(t),  \ y \in \Rset^  ,
\end{equation}
for non-decreasing $r$.  Similarly we have that
\begin{equation}
\Phi^{h,r}_{(a,b]}(y)= \int_{(a,b]} \Phi(-h(t)y)|dr(t)|, \ \  y\in \Rset^d ,
\end{equation}
for non-increasing $r$, because for $0<u<w$, we have $\mathcal{L}(Y_{\nu}(u)-Y_{\nu}(w))=
(\nu^{-})^{\ast(w-u)}$ where $\nu^{-}:=\mathcal{L}(-Y_{\nu}(1))$. Consequently we have distributional equality of two processes:  $(-Y_{\nu}(t), t\ge0)\stackrel{d}{=}(Y_{\nu^-}(t),
t\ge0)$.

\medskip
2 d). \underline{\it{Proofs of Theorem 1 , Corollaries 1 and 2.}}

\ For $\nu\in\mathcal{D}^{h,r}_{(a,b]}$ and its L\'evy
exponent $\Phi$ let us define the script mapping
$\mathcal{I}^{h,r}_{(a,b]}$ as follows
\begin{equation}
\mathcal{I}^{h,r}_{(a,b]}(\Phi)(y):= \Phi^{h,r}_{(a,b]}=
\int_{(a,b]}\Phi(\pm h(s)y)d(\pm)r(s),
\end{equation}
where the sign minus is in the case of decreasing time change $r$.
Then to justify (2) it is enough to notice that
\begin{multline}
\mathcal{I}_{(a_1, b_1]}^{h_1, \rho_1}(\mathcal{I}_{(a_2,
b_2]}^{h_2, \rho_2}(...(\mathcal{I}_{(a_m,
b_m]}^{h_m,\rho_m}(\Phi))...))(y)\\ =
\int_{(a_1,b_1]}\int_{(a_1,b_2]}...\int_{(a_m,b_m]}\Phi\big(h_1(t_1)\,h_2(t_2)\,...\,h_m(t_m)\,y)\big)\,dr_m(t_m)...dr_2(t_2)dr_1(t_1)\\
=\int_{(\textbf{a},\textbf{b}]}\Phi\big(h_1 \otimes...\otimes
h_m(s)\,y\big) \boldsymbol{\rho}(ds) =
\int_{(c,d]}\Phi(t\,y)(\textbf{h} \boldsymbol{\rho})(dt),
\end{multline}
which follows from the Fubini and the image measure theorems.

To conclude the second  part of Theorem 1 (commutativity) one needs to note that
\[
h_1 \otimes...\otimes
h_m\,(\rho_1\times\rho_2\times...\times\rho_m)=h_{\sigma(1)}
\otimes...\otimes
h_{\sigma(m)}\,(\rho_{\sigma(1)}\times\rho_{\sigma(2)}\times...\times\rho_{\sigma(m)}),
\]
for any permutation $\sigma$ of $1,2,..., m$.

Corollary 1  follows from the definition of tensor product of functions.
\begin{rem}
Corollary 1, for $h_i(t)=|t|$ and standard normal random variables $Z_i$ was investigated by Aoyama (2009) via polar decomposition of L\'evy spectral measures.
\end{rem}

For a proof of Corollary 2 note  that  the equality
\[
I^{h,r}_{(a,b]}(\mathcal{D}^{h,r}_{(a,b]})=I^{h_1,r_1}_{(a_1,b_1]}\big(
I^{h_2,r_2}_{(a_2,b_2]}(\mathcal{D}^{h,r}_{(a,b]})\big) \ \mbox{ implies}\ \ I^{h_2,r_2}_{(a_2,b_2]}(\mathcal{D}^{h,r}_{(a,b]})\subset
\mathcal{D}^{h_1,r_1}_{(a_1,b_1]}.
\]
Hence  $I^{h,r}_{(a,b]}(\mathcal{D}^{h,r}_{(a,b]})\subset
I^{h_1,r_1}_{(a_1,b_1]}(\mathcal{D}^{h_1,r_1}_{(a_1,b_1]})$. By the commutative    property
we also  get
$I^{h,r}_{(a,b]}(\mathcal{D}^{h,r}_{(a,b]})\subset
I^{h_2,r_2}_{(a_2,b_2]}(\mathcal{D}^{h_2,r_2}_{(a_2,b_2]})$, which
completes a proof of Corollary 2.

\medskip
2 e). \underline{\it{Proof of Theorem 2.}}

\medskip
 Part (a). The homomorphism property of $I^{h,r}_{(a,b]}$, that is, the equality
\[
I^{h,\, r}_{(a,b]}\,(\nu_1\ast\nu_2))=I^{h,\,
r}_{(a,b]}\,(\nu_1)\ast I^{h,\, r}_{(a,b]}\,(\nu_2),
\]
in terms of the corresponding L\'evy exponents, follows from (8) or (9).

For the continuity, let us note that $0\leq |r(b)-r(a)|<\infty$  and
the cadlag property imply that functions $t\to Y(r(t))$ are bounded
and with at most countable many discontinuities; cf. Billingsley
(1968), Chapter 3, Lemma 1. Furthermore, the mapping
\begin{multline} D_{\Rset^d}[a,b] \ni y\to
\int_{(a,b]}h(t)dy(r(t)):=
\\ h(b)y(r(b))-h(a)y(r(a))- \int_{(a,b]}y(r(t)-)dh(t) \in \Rset^d,
\end{multline}
is continuous in Skorohod topology (for details see Billingsley
(1968), p. 121.). Furthermore, if $\nu_n \Rightarrow \nu$ then
$(Y_{\nu_{n}}(t), a \le t \le b) \Rightarrow (Y_{\nu}(t), 0\le t\le
b)$ in  Skorohod space $D_{\Rset^d}[a,b]$ of cadlag functions.  Consequently, we have
\[
\mathcal{L}\Big(\int_{(a,b]}h(t)dY_{\nu_{n}}(r(t))\Big)\Rightarrow\mathcal{L}\Big(\int_{(a,b]}h(t)dY_{\nu}(r(t))\Big),
\]
which proves the continuity of mappings  $I^{h,r}_{(a,b]}$ and completes the proof of  part (a) of Theorem 2.

Equality  in parts (b), (c) and (d) follow from (8) and (9).

\medskip
\medskip
\textbf{3. Applications and illustrations of Theorem 1.}

\medskip
We begin with the following auxiliary fact.
\begin{lem}
Let $ h_1(t):= e^{-t},  r_1(t):=t, h_2(s):=s$ and
$r_2(s):=1-e^{-s}$, \qquad  \ $0<s, t<\infty$. Then the corresponding measures are:
$d\rho_1(t)=dt, \, d\rho_2(s)=e^{-s}ds$ and
$d\boldsymbol{\rho}(t,s)= d(\rho_1\times
\rho_2)(t,s)=e^{-s}\,dt\,ds$. Finally, for the image measure
$\textbf{h}\,\,\boldsymbol{\rho}(dw)=(h_1\otimes h_2) (\rho_1\times
\rho_2)(dw)=\frac{e^{-w}}{w}dw$.
\end{lem}
\emph{Proof.} For Borel measurable, bounded and non-negative
functions $g$ we have
\begin{multline*}
\int_0^{\infty}g(u)(h_1\otimes h_2) (\rho_1\times
\rho_2)(du)=\int_0^{\infty}\int_0^{\infty}g((h_1\otimes
h_2)(t,s))\rho_1(dt)\rho_2(ds)\\
=\int_0^{\infty}\int_0^{\infty}g(e^{-t}\,s)dt\,e^{-s}ds=
\int_0^{\infty}(\int_0^{s}g(w)\frac{1}{w}\,dw)\,e^{-s}ds=
\int_0^{\infty}g(s)\,\frac{e^{-s}}{s}\,ds,
\end{multline*}
which completes the proof of Lemma 1.

From Theorem 1 and  Lemma 1 we conclude the following.
\begin{ex} For $\nu\in ID_{\log}$  we have
\[
 I^{t,\,1-e^{-t}}_{(0,\infty)}\big(I^{e^{-s},\, s}_{(0,\infty)}(\nu)\big)=
I^{e^{-s},\,s}_{(0,\infty)}\big(I^{t,\,1-e^{-t}}_{(0,\infty)}(\nu)\big)=I^{-w,\,\Gamma(0;w)}_{(0,\infty)}(\nu)=I^{w,\,\Gamma(0;w)}_{(0,\infty)}(\nu^-).
\]
Moreover, $\Gamma(0; w)= (h_1\otimes h_2) (\rho_1\times
\rho_2)(\{x:x>w\})=\int_w^{\infty}\frac{e^{-s}}{s}\,ds$, for $w>0$.
\end{ex}

\begin{rem} (a) For the Euler constant \textbf{C} we have
\[
- \Gamma(0;w)=Ei(-w)= \textbf{C}+ ln\,w
+\int_0^w\frac{e^{-t}-1}{t}dt, \ \mbox{for} \  w>0,
\]
where $Ei$ is the special \emph{exponential-integral} function; cf.
Gradshteyn-Ryzhik (1994), formulas  8.211 and 8.212.

(b) Recall that the class 
$I^{t,\,1-e^{-t}}_{(0,\infty)}(ID)\equiv \mathcal{E}$ was introduced
in Jurek (2007), where the mapping $I^{t,\,1-e^{-t}}_{(0,\infty)}$
was denoted by $\mathcal{K}^{(e)}$; \  (here $ (e)$ stands  for the exponential 
distribution). 

More importantly, the class $\mathcal{E}$
was related to the class of Voiculescu $\boxplus$ free-infinitely
divisible measures; cf. Corollary 6 in Jurek (2007). Note also that
$I^{t,\,1-e^{-t}}_{(0,\infty)}=I^{-\log s,\,s}_{(0,1]}$ and thus it
coincides with the upsilon mapping $\Upsilon$ studied in
Barndorff-Nielsen, Maejima and Sato (2006).

(c) Similarly $I^{e^{-s},\,s}_{(0,\infty)}(ID_{\log})\equiv L$
coincides with the L\'evy class of selfdecomposable probability
measures; cf. Jurek-Vervaat (1983), Theorem 3.2 or 
Jurek-Mason (1993), Theorem 3.6.6.

(d) Finally we get identity
$I^{e^{-s},\,s}_{(0,\infty)}\big(I^{t,\,1-e^{-t}}_{(0,\infty)}(ID_{\log})\big)
\equiv T$, which is the Thorin class; cf. Grigelionis (2007),
Maejima and Sato (2009) or Jurek (2011).
\end{rem}

From Corollary 2 ,  Remark 3 (c) and (d) we infer the following inclussion.
\begin{cor}
For the three classes: Thorin class T, L\'evy class L
(selfdecomposable measures) and $\mathcal{E}$  of probability measures on $\Rset^d$, we have that $T\subset
L\cap\mathcal{E}$.
\end{cor}
This inclusion was first noticed in
Barndorff-Nielsen, Maejima and Sato (2006) and also in Remark 2.3 in
Maejima-Sato (2009) but by using different methods. See note  (c)  in Concluding Remarks.

Finally, we give additionally three examples of compositions of random integral mappings.
\begin{ex} For $\beta>0$ we have
\begin{equation}
I^{t^{1/\beta},\, t}_{(0,1]}\circ I^{s^{1/2\beta},\,
s}_{(0,1]}=I^{w,\, 2w^{\beta}(1-(1/2)
w^{\beta})}_{(0,1]}=I^{(1-\sqrt{t})^{1/\beta},\, t}_{(0,1]}\,.
\end{equation}
Or equivalently, for Lebesque measure $l_1$ on the unit interval and
$ 0< w \le 1$ we get
\begin{multline*}
(t^{1/\beta}\otimes s^{1/(2\beta)})(l_1 \times l_1)(dw) \\ = id
^{\otimes 2}(\beta t^{\beta-1}dt \times 2 \beta s^{2
\beta-1}dt)(dw)=2 \beta w^{\beta-1}(1-w^{\beta})\,dw.
\end{multline*}
\end{ex}
\emph{Proof.} As in Example 1, it simply follows from Theorem 1 and
identity (8),  because all time change functions are strictly
increasing on the unit interval.
\begin{ex}
For $\beta >0$
\[
I^{t^{1/\beta},\, t}_{(0,1]}\circ I^{e^{-s},\,s}_{(0,\infty)}=
I^{e^{-s},\,\,s+\beta^{-1}e^{-\beta
s}-\beta^{-1}}_{(0,\infty)}=I^{-w,\,\, \beta^{-1}w^{\beta}-\log
w-\be^{-1}}_{(0,1]}.
\]
Or equivalently, for $0<w\le1$
\begin{equation*}
(t^{1/\beta}\otimes e^{-s})(l_1 \times l)(dw)=
(\beta^{-1}w^{\beta}-\log w-\be^{-1})dw.
\end{equation*}
\end{ex}
\noindent This is a consequence of Theorem 1. Also cf. Czyżewska-Jankowska
and Jurek (2011), Proposition 2.
\begin{ex} For $\al\in \Rset$, let $\Gamma(\alpha; x):=\int_x^\infty t^{\alpha-1}e^{-t}dt, x>0$ be the incomplete Euler  function. Then we have
\[
I^{t,\, \Gamma(\al; t)}_{(0,\infty)}\circ
I^{e^{-s},\,s}_{(0,\infty)}=
I^{t,\,\int_t^{\infty}s^{-1}\Gamma(\al;s)ds}_{(0,\infty)},
\]
\end{ex}
\noindent which follows from Theorem 1.

\medskip
\textbf{4. Concluding Remarks.} 

(a) Let $E$ be a real separable Banach space with the dual $E^\prime$ and a bilinear form $ <\cdot ,\cdot>: E^\prime \times  E \to \Rset$,  (a scalar product when $E$ is a Hilbert space). Then we define a convolution $\ast$ of measures $\mu$,  a characteristic functional (Fourier transfrom) $\hat{\mu}$  and, in particular, the notion of infinite divisibility of (Borel) probability measures in $E$ and the  convolution semigroup $ID(E)$. The L\'evy-Khintchine representation (4)  holds true with $\Rset^d$ replaced by $E$. However, for infinite dimesional Banach spaces,  the condition $\int_E \min(||x||^2, 1)M(dx)<\infty$ is neither sufficient or necessary for $M$ to be  a  L\'evy (spectral) measure.  Moreover, in general case of E,  convergence of characteristic functionals $\hat{\mu_n}(y)\to\hat{\mu}(y), y \in E^\prime, $ does not imply the weak convergence $\mu_n \Rightarrow \mu$.
cf. Araujo-Gine (1980) or Parthasarathy (1967).   

(b) Because in this note we did not use any of those two exceptions our results are valid for measures on real  separable Banach spaces E;  (and E-valued L\'evy processes $Y$). Note that in the proof of continuity of random integral mappings (Theorem 2) we did not use characteristic functional arguments.

(c) Corollary 3,  inclusions of three ranges of a given three random integral mappings, is valid for measures on Banach spaces. However, we do not define Thorin class $T$ via properties of  L\'evy measures  M; cf. for instance  Jurek (2011) or Grigelionis (2007).

\medskip
\textbf{Acknowledgements.} Comments of a Reviewer helped to reorganized this note. In particular, the discussion of infinite dimensional case was postponed to the last section of this paper.

\medskip
\begin{center}
\textbf{References}
\end{center}
\noindent [1] T. Aoyama (2009), Nested subclasses of the class of type G
selfdecomposable distributions on $\Rset^d$, \emph{Probab. Math.
Stat.}, vol. 29,  pp. 135-154.

\noindent  [2]  A. Araujo, and E. Gin\'e (1980), \textit{The Central Limit Theorem for Real
and Banach Valued Random Variables}, Wiley, New York.

\noindent [3] O. E. Barndorff-Nielsen, M. Maejima and K. Sato (2006),
Some classes of multivariate infinitely divisible distributions
admitting integral representations, \emph{Bernoulli} \textbf{12},
pp. 1-33.

\noindent[4] P. Billingsley (1968), \emph{Convergence of probability
measures}, Wiley, New York.

\noindent [5]  A. Czyżewska-Jankowska and Z. J. Jurek (2011),
Factorization property of generalized s-selfdecomposable measures
and class $L^f$ distributions, \emph{Theory Probab. Appl.} vol.
\textbf{55}, no 4, pp. 692-698.

\noindent[6] I. S. Gradshteyn and I. M. Ryzhik (1965).
\emph{Tables of integrals, series, and products}, Academic Press,
New York.

\noindent[7]  B. Grigelionis (2007), Extended Thorin classes and
stochastic integrals, \emph{Liet. Matem. Rink.} \textbf{47} , pp.
497 -- 503.

\noindent [8]  Z.J. Jurek (2007), Random integral representations for free-infinitely divisible
and tempered stable distributions, \emph{Stat.\& Probab. Letters},
\textbf{77} no. 4, pp. 417-425.

\noindent[9]  Z. J. Jurek (2011), The Random Integral Representation Conjecture:  a quarter of a century later,
\emph{Lithuanian Math. Journal}, \textbf{51}, no 3, 2011, pp.
362-369.

\noindent[10] Z. J. Jurek (2014), Remarks on the factorization property of some random integrals,  \emph{Stat.\& Probab. Letters},
\textbf{94}, pp. 192-195.

\noindent [11]  Z. J. Jurek and J.D. Mason (1993), \textit{Operator-limit Distributions
in Probability Theory}, Wiley Series in Probability and Mathematical
Statistics, New York.

\noindent[12]   Z.J. Jurek and W. Vervaat (1983), An integral representation for
selfdecomposable Banach space valued random variables, \textit{Z.
Wahrsch. verw. Gebiete} \textbf{62}, 247--262.

\noindent [13]  M. Maejima, V.Perez Abreu, K.I. Sato (2012), A class of
multivariate infinitely divisble distributions related to arcsine
density, \emph{Bernoulli}, vol. 18 no. 2, pp. 476-495.

\noindent [14] M. Meerschaert and P. Scheffler (2001), \textit{Limit distributions for sums of independent random vectors}, J. Wiley $\&$ Sons,  New York.

\noindent [15]  K. R. Parthasarathy (1968), \textit{Probability measures on
metric spaces}, Academic Press, New York and London, 1968.

\noindent [16]  K. Sato (2006), Two families of improper stochastic
integrals with respect to L\'evy processes, \emph{ALEA, Lat. Am. J.
Probab. Math. Stat.} vol. 1, pp. 47-87.

\end{document}